\documentclass[12pt,twoside,sort&compress ]{article} % say
\usepackage{mathrsfs,amsfonts,amsmath,amsbsy,tikz}\usepackage{fontenc}\usepackage{textcomp}
\usepackage{graphicx}
\usepackage{amssymb}
\usepackage{graphicx}

\def\ge{\geq}

\def\bbb{\begin{eqnarray*}}

\def\eee{\end{eqnarray*}}

\begin{document} % say
\baselineskip=18pt
\begin{center}
\vspace{-0.6in}{\large \bf Weakly mixing, topologically weakly mixing, and sensitivity \\[0.1in]
 for non-autonomous discrete systems}\\ [0.2in]
 Hua Shao,\ \ Yuming Shi$^{*}$, \ \ Hao Zhu\\
\vspace{0.15in}  Department of Mathematics, Shandong University\\
Jinan, Shandong 250100, P. R. China\\

\footnote{$^{*}$ The corresponding author.}
\footnote{ Email addresses: huashaosdu@163.com (H. Shao), ymshi@sdu.edu.cn (Y. Shi), haozhu@mail.sdu.edu.cn (H. Zhu).} \

\end{center}

{\bf Abstract.} This paper is concerned with relationships of weakly mixing, topologically weakly mixing, and sensitivity
for non-autonomous discrete systems. It is shown that weakly mixing implies topologically weakly mixing and sensitivity for
measurable systems with a fully supported measure; and topological weakly mixing implies sensitivity for general dynamical
systems. However, the inverse conclusions are not true and some counterexamples are given. The related existing results for autonomous discrete
systems are generalized to non-autonomous discrete systems and their conditions are weaken.
\medskip

\medskip

{2000 {\bf \it Mathematics Subject Classification}}: 37B55; 37A25; 49K40.

\bigskip

\noindent{\bf 1. Introduction}\medskip

Chaos is a universal dynamical behavior of nonlinear dynamical systems and one of the central topics of research on nonlinear
science. It is well known that sensitivity characterizes the unpredictability of chaotic phenomena,
and is the essential condition of various definitions of chaos. Therefore, the study on sensitivity has attracted a lot of attention
from many scholars [1, 2, 5, 7, 8, 10, 12].

In 2002, Abraham et al. proved that if a measure-preserving map $f$ on a metric probability space $X$ with a fully supported measure is
either topologically mixing or weakly mixing, and satisfies that for any nonempty open set $U\subset X$, there is a sequence $\{n_k\}_{k=0}^{\infty}$
with positive upper density such that $U\cap(\cap_{k\geq0} f^{-n_k}U)\neq\emptyset$, then $f$ is sensitive [1]. In 2004, He et al. relaxed the conditions of the above result and showed that if a measure-preserving map (resp. a measure-preserving semi-flow) on a metric probability space with a fully supported measure is weakly mixing, then it is sensitive [8]. In 2010, Li and the second author of the present paper gave several sufficient conditions of sensitivity for maps and semi-flows, and showed that if a measure-preserving map (resp. a measure-preserving semi-flow) on a metric probability space with a fully supported measure is topologically strongly ergodic, then it is sensitive [12].

Since many complex systems occurring in the real-world problems such as physical, biological, and economical problems are necessarily described by non-autonomous discrete systems, which are generated by iteration of a sequence of maps in a certain order, many scientists and mathematicians focused on complexity of non-autonomous discrete systems recently [3, 6, 10, 11, 14--18, 20]. In 1996, Kolyada and Snoha introduced the concept of topological entropy for a non-autonomous discrete system and studied its properties [11]. In 2006, Tian and Chen extended the concept of chaos in the sense of Devaney to non-autonomous discrete systems [18]. In 2009, related concepts of chaos, such as topological transitivity, sensitivity, chaos in the sense of Li-Yorke, Wiggins, and Devaney were extended to general non-autonomous discrete systems [15]. Recently, some sufficient conditions of sensitivity for general non-autonomous discrete systems are presented in [10]. Motivated by the above results, we shall try to investigate the relationships of weakly mixing, topologically weakly mixing, and sensitivity for non-autonomous discrete systems in the present paper.

The rest of the paper is organized as follows. In Section 2, some basic concepts and useful lemmas are presented.
In Section 3, relationships of weakly mixing, topologically weakly mixing, and sensitivity for non-autonomous discrete systems are discussed.
\bigskip

\noindent{\bf 2. Preliminaries }\medskip

In this section, some basic concepts for non-autonomous discrete systems are introduced,
including topological transitivity, topologically weakly mixing, topologically mixing, weakly mixing, and sensitivity.
In addition, some useful lemmas are also presented.\medskip

Let $\mathbb{N}$ and $\mathbb{Z^{+}}$ denote the set of all nonnegative integers and that of all positive integers, respectively, and set $N_n:=\{0,1,\cdots, n-1\}$. Define $\mathbb{R}^{+}:=[0,+\infty)$ and $N(U,V):=\{n\in \mathbb{Z^{+}}: f_0^{n}(U)\cap V\neq\emptyset\}$
for any two nonempty open sets $U,V$.\medskip

We shall consider the following non-autonomous discrete system in the present paper:
\vspace{-0.5cm}$$x_{n+1}=f_n(x_n),\;\;\;\; n\geq 0, \eqno(2.1)\vspace{-0.2cm}$$
where $f_n: X\to X$ is a map for each $n\geq 0$ and $(X,d)$ is a metric space.\medskip

For any fixed $x_0\in X$, $\{x_n\}_{n=0}^{\infty}$ is called the (positive) orbit of system (2.1) starting from $x_0$
and $x_n=f_0^n(x_0),\;n\geq0$, where $f_0^n:=f_{n-1}\circ\cdots\circ f_{0}$.
Furthermore, by $B_{\epsilon}(x_0)$ denote the open
ball in $X$ of radius $\epsilon$, centered at $x_0$.\medskip

\noindent{\bf Definition 2.1} [15, Definition 2.2]. Let $A$ be a nonempty subset of $X$. System (2.1) is said to be topologically
transitive in $A$ if for any two nonempty relatively open subsets $U$ and $V$ with respect to $A$, $N(U,V)\neq\emptyset$.\medskip

\noindent{\bf Definition 2.2.} Let $A$ be a nonempty subset of $X$. System (2.1) is said to be topologically weakly mixing
in $A$ if for any four nonempty relatively open subsets $U_1, V_1, U_2, V_2$ with respect to $A$,
$N(U_1,V_1)\cap N(U_2,V_2)\neq\emptyset$.\medskip

\noindent{\bf Definition 2.3} [14, Definition 2.3]. Let $A$ be a nonempty subset of $X$. System (2.1) is said to be topologically
mixing in $A$ if for any two nonempty relatively open subsets $U$ and $V$ with respect to $A$, there is a positive integer $N$ such that
$N(U,V)\supset[N, +\infty)\cap\mathbb{Z^{+}}$.\medskip

Note that topologically mixing implies topologically weakly mixing, and topologically weakly mixing implies topological
transitivity for system (2.1) in $A$.\medskip

\noindent{\bf Definition 2.4} [15, Definition 2.3]. Let $A$ be a nonempty subset of $X$. System (2.1) is said to have sensitive
dependence on initial conditions in $A$ if there exists a constant $\delta_0> 0$  such that for any
$x_0\in A$ and any neighborhood $U$ of $x_0$, there exist $y_0\in A\cap U$ and a positive integer $N$
such that $d(f_0^N(x_0), f_0^N(y_0))>\delta_0$. The constant $\delta_0$ is called a sensitivity constant of system (2.1)
in $A$.\medskip

Let $\mathcal{B}(X)$ be the $\sigma$-algebra of Borel subsets of $X$, and $\mu$ be a finite measure
of measurable space $(X,\mathcal{B}(X))$.\medskip

\noindent{\bf Definition 2.5} [10, Definition 2.5]. System (2.1) is measurable if $f_n$ is measurable
on $(X,\mathcal{B}(X),\mu)$ for each $n\geq0$.\medskip

\noindent{\bf Definition 2.6.} Let system (2.1) be measurable on $(X ,\mathcal{B}(X),\mu)$. It is said to be weakly mixing
if for any $A$, $B\in \mathcal{B}(X)$,
\vspace{-0.2cm}$$ \lim\limits_{n\to\infty}\frac{1}{n}\sum_{i=0}^{n-1}\big|\mu\big(A\cap f_0^{-i}(B)\big)-\mu(A)\mu(B)\big|=0.\eqno(2.2)\vspace{-0.2cm}$$

Let $S\subset\mathbb{N}$ and $|S|$ be the cardinality of $S$. Set
\vspace{-0.2cm}$$\overline{d_{0}}(S):=\limsup\limits_{n\rightarrow\infty}\frac{1}{n}|S\cap N_{n}|,\;\;\underline{d_{0}}(S):=\liminf\limits_{n\rightarrow\infty}\frac{1}{n}|S\cap N_{n}|.
\vspace{-0.2cm}$$
$\overline{d_{0}}(S)$ and $\underline{d_{0}}(S)$ are called the upper density and the lower density of $S$, respectively.
If $\overline{d_{0}}(S)=\underline{d_{0}}(S)=:d_{0}(S)$, then $d_{0}(S)$ is called the density of $S$.\medskip

\noindent{\bf Lemma 2.1} [13, Lemma 2.6.2]. {\it Let $\{a_i\}_{i=0}^{\infty}$ be a bounded sequence of nonnegative numbers.
Then $\lim_{n\rightarrow\infty}\frac{1}{n}\sum_{i=0}^{n-1}a_i=0$
if and only if there exists a subset $E\subset \mathbb{N}$ of density zero such that
$\lim\limits_{\scriptstyle n\rightarrow\infty \atop \scriptstyle n\overline{\in}E}a_n=0$.}\medskip

\noindent{\bf Lemma 2.2.} {\it Let system {\rm (2.1)} be measurable on $(X ,\mathcal{B}(X), \mu)$.
For any $A,B\in\mathcal{B}(X)$, {\rm(2.2)} holds if and only if there exists $J\subset \mathbb{N}$ with $d_{0}(J)=1$ such that
\vspace{-0.2cm}$$\lim\limits_{\scriptstyle n\rightarrow\infty \atop \scriptstyle n\in J}\mu\big(A\cap f_0^{-n}(B)\big)=\mu(A)\mu(B).\vspace{-0.7cm}$$}\medskip

\noindent{\bf Proof.} It can be directly derived by Lemma 2.1.
\bigskip

\noindent{\bf 3. Relationships of weakly mixing, topologically weakly mixing, and sensitivity}\medskip

In this section, we shall investigate some relationships of weakly mixing, topologically weakly mixing,
and sensitivity for system (2.1).\medskip

\noindent{\bf Theorem 3.1.} {\it If system {\rm (2.1)} is topologically weakly mixing, then it is sensitive in $X$.}\medskip

\noindent{\bf Proof.} Fix any two different points $x_0, y_0\in X$. Denote $d(x_0,y_0):=8\delta>0$. Then for any $x\in X$,
\vspace{-0.2cm}$$d(x_0,x)\geq 4\delta\;\;{\rm or}\;\; d(y_0,x)\geq 4\delta.\vspace{-0.2cm}$$
Assume that $d(x_0,x)\geq 4\delta$. Since system (2.1) is topologically weakly mixing in $X$,
for any $0<\epsilon<\delta$, there exists a positive integer $n_1$ such that
\vspace{-0.2cm}$$f_0^{n_1}(B_{\epsilon}(x))\cap B_{\delta}(x_0)\neq \emptyset, \;\;f_0^{n_1} (B_{\epsilon}(x)) \cap B_{\delta}(x)\neq \emptyset.\vspace{-0.2cm}$$
So there exist two points $\hat x,\hat y\in B_{\epsilon}(x)$ such that $f_0^{n_1}(\hat x) \in B_{\delta}(x_0)$ and
$f_0^{n_1}(\hat y)\in B_{\delta}(x)$. Therefore,
\vspace{-0.2cm}$$d(f_0^{n_1}(\hat x),f_0^{n_1}(\hat y))\geq d(x_0,x)-d(x,f_0^{n_1}(\hat y))-d(x_0,f_0^{n_1}(\hat x))\geq 4\delta-\delta-\delta=2\delta.\vspace{-0.2cm}$$
Hence, $d(f_0^{n_1}(\hat x),f_0^{n_1}(x))\geq \delta$ or $d(f_0^{n_1}(\hat y),f_0^{n_1}(x))\geq \delta$.
Therefore, system (2.1) is sensitive in $X$. This completes the proof.\medskip

\noindent{\bf Remark 3.1.} Note that topologically mixing implies topologically weakly mixing for system (2.1). So, topologically mixing
implies sensitivity for system (2.1) by Theorem 3.1. This implies that condition (i) in the definition of chaos
in the strong sense of Wiggins in [14, Definition 2.5] is redundant.\medskip

\noindent{\bf Theorem 3.2.} {\it Let system {\rm (2.1)} be measurable on $(X ,\mathcal{B}(X ), \mu)$  with $supp\,\mu=X$.
If {\rm(2.2)} holds for any two nonempty open subsets $A$ and $B$ of $X$, then
\begin{itemize}\vspace{-0.2cm}
\item[{\rm (i)}] system {\rm(2.1)} is topologically weakly mixing;
\item[{\rm (ii)}] system {\rm(2.1)} is sensitive.

\end{itemize}\vspace{-0.2cm}
}\medskip

\noindent{\bf Proof.} First, we show that assertion (i) holds.  For any nonempty open subsets $A_1$, $B_1$, $A_2$, $B_2$ of $X$, by Lemma 2.2
there exist two sets $J_1, J_2 \subset \mathbb{N}$ with $d_{0}(J_1)=d_{0}(J_2)=1$ such that
\vspace{-0.2cm}$$\lim\limits_{\scriptstyle n\rightarrow\infty \atop \scriptstyle n\in J_i} \mu(A_i\cap f_0^{-n}B_i)=\mu(A_i)\mu(B_i),\;i=1,2.\vspace{-0.2cm}$$
Since $supp\,\mu=X$, $\mu(A_i)\mu(B_i)>0$, $i=1,2$. So, there exists an integer $N>0$ such that
for each $n>N$ with $n\in J_i$, $\mu(A_i\cap f_0^{-n}B_i)>0$ and thus $f_0^{n}(A_i)\cap B_i\neq\emptyset$, $i=1,2$.

Note that $d_{0}(J_i\backslash\{0,\cdots, N\})=d_{0}(J_i)=1,\;i=1,2$. Thus,
$(J_1\cap J_2)\backslash\{0,\cdots, N\}\neq\emptyset.$
In fact, suppose that $(J_1\cap J_2)\backslash\{0,\cdots, N\}=\emptyset$.
Set $\hat{J}_1:=J_1\backslash\{0,\cdots, N\},\;
\hat{J}_2:=J_2\backslash\{0,\cdots, N\}.$
Then $\hat{J}_1\cap\hat{J}_2=\emptyset$ and $d_{0}(\hat{J}_1)=d_{0}(\hat{J}_2)=1$. Thus, one has that
\vspace{-0.2cm}$$\begin{array}{llll}\lim\limits_{n\to\infty}\frac{1}{n}|(\hat{J}_{1}\cap N_{n})\cup(\hat{J}_{2}\cap N_{n})|
&=&\lim\limits_{n\to\infty}\frac{1}{n}|\hat{J}_1\cap N_{n}|+\lim\limits_{n\to\infty}\frac{1}{n}|\hat{J}_2\cap N_{n}|\\[1.5ex]
&=&d_{0}(\hat{J}_1)+d_{0}(\hat{J}_2)=2.\end{array}\vspace{-0.2cm}$$
On the other hand, one gets that
\vspace{-0.2cm}$$\lim_{n\to\infty}\frac{1}{n}|(\hat{J}_1\cap N_{n})\cup(\hat{J}_2\cap N_{n})|\leq\lim_{n\to\infty}\frac{1}{n}|N_{n}|=1,\vspace{-0.2cm}$$
which is a contradiction. So, there exists $n_0\in(J_1\cap J_2)\backslash\{0,\cdots, N\}$ and thus $n_0\in N(A_1,B_1)\cap N(A_2,B_2)\neq\emptyset$.
Therefore, system {\rm(2.1)} is topologically weakly mixing.

Assertion (ii) holds by assertion (i) and Theorem 3.1. This completes the proof.\medskip

\noindent{\bf Remark 3.2.} Assertion (ii) in Theorem 3.2 can also be proved in a direct way using Theorem 3.1 in [10].
\medskip

\noindent{\bf Theorem 3.3.} {\it Let system {\rm (2.1)} be measurable on $(X,\mathcal{B}(X),\mu)$ with $supp\,\mu=X$.
If system {\rm (2.1)} is weakly mixing, then
\begin{itemize}\vspace{-0.2cm}
\item[{\rm (i)}] system {\rm(2.1)} is topologically weakly mixing;
\item[{\rm (ii)}] system {\rm(2.1)} is sensitive.
\end{itemize}\vspace{-0.2cm}
}\medskip

\noindent{\bf Proof.} Since system {\rm (2.1)} is weakly mixing, {\rm(2.2)} holds for any two nonempty open subsets $A$ and $B$ of $X$, and thus
the two assertions hold by Theorem 3.2. This completes the proof.\medskip

\noindent{\bf Remark 3.3.} Assertion (ii) in Theorems 3.2 and 3.3 extends that of [8, Theorem A] for autonomous discrete systems to
non-autonomous discrete systems. The condition that the map is measure-preserving in [8, Theorem A] is removed in Theorems 3.2 and 3.3;
and the condition of weakly mixing in [8, Theorem A] is relaxed since it is only required that {\rm(2.2)} holds
for any two nonempty open subsets $A$ and $B$ of $X$ in Theorem 3.2.\medskip

By the results of Theorems 3.1--3.3, one gets that weakly mixing implies topologically weakly mixing
for measurable systems with a fully supported measure, and topologically weakly mixing implies sensitivity for general dynamical systems.
It is natural to ask whether their reverse conclusions hold. That is, does topologically weakly mixing imply weakly mixing? and does sensitivity imply
topologically weakly mixing? For the first question, one counterexample has been given, see [4, p. 32] and [19]. For the second
question, a counterexample is given as follows.
\medskip

\noindent{\bf Example 3.1.} Let $X=[0,3/2]$, $f_n=f$ for each $n\geq0$, and
\vspace{-0.1cm}$$\vspace{-0.2cm}
f(x) =\begin{cases}
2x & \text{ if } x\in [0,1/2], \\
                     2(1-x) & \text{ if } x\in(1/2,1],\\
                     2(x-1) & \text{ if } x\in(1,3/2].
\end{cases}$$\vspace{-0.1cm}

Clearly, $f$ maps $X$ into itself.
Since $f^{n}((0,1))\cap(1,3/2)=\emptyset$ for each $n\ge 1$,
$f$ is not topologically transitive in $X$ and thus it is not topologically weakly mixing in $X$.
On the other hand, one can easily verify that $f$ is sensitive in $X$ since the tent map is sensitive in [0, 1],
and $f(1,3/2]\subset[0, 1]$.
Therefore, $f$ is sensitive in $X$.
\bigskip

\noindent{\bf Acknowledgment}\medskip

This research was supported by the NNSF of China (Grant 11571202).

\end{document}